\documentstyle[12pt]{amsart}

\newtheorem{thm}{Theorem}

\newcommand{\C}{\Bbb C}
\newcommand{\CP}{{\Bbb C}P}

\newcommand{\OO}{\cal O}

\newcommand{\Z}{\Bbb Z}
\newcommand{\cp}{\C P^2}

\begin{document}

{\hfill NI94036}

\bigskip\bigskip

\title{Non--trivial harmonic spinors on\\ generic algebraic surfaces}

\author{D. Kotschick}
\address{Mathematisches Institut\\ Universit\"at Basel\\
Rheinsprung 21\\ 4051 Basel\\ Switzerland}
\thanks{This note was written while the author was an
EPSRC Visiting Fellow at the Isaac Newton Institute for
Mathematical Sciences in Cambridge}
\email{dieter@@math.unibas.ch}

\maketitle

\bigskip\bigskip

\centerline{({\it to appear in Proc. Amer. Math. Soc.})}

\bigskip\bigskip

For every closed Riemannian spin manifold its Dirac operator is a
selfadjoint linear elliptic operator acting on sections of the
spinor bundle. Hitchin~\cite{hi} proved that the dimension of
its kernel, the space of harmonic spinors,
only depends on the conformal class of the Riemannian metric,
and that it varies with the conformal class. However, the
variation is not understood.

Investigating this variation is particularly interesting
for manifolds of dimension $4k$. In that case the Dirac operator
interchanges spinors of different chirality and the index
of the (half--) Dirac operator $D\colon V_+ \rightarrow V_-$
given by the $\hat{A}$--genus can be non--zero. We say that a
conformal structure admits non--trivial harmonic spinors if
the dimension of the space of harmonic spinors is larger than
$\vert index(D)\vert$. This is equivalent to the existence
of harmonic spinors of both chiralities. The natural problem
to consider is:

\medskip\noindent
{\bf Problem 1.} {\sl Does there always exist a conformal
structure for which there are no non--trivial harmonic
spinors? Are such conformal structures generic?}
\medskip

For manifolds of dimension $4$, this question
is related to the applications of Witten's
monopole equations \cite{wi}.

In this note we address the complex analytic analog of
Problem 1. For a complex manifold $X$, the existence of a
spin structure is equivalent to the existence of a
half--canonical divisor $\frac{1}{2}K\in H^2 (X,\Z )$.
Hodge theory gives canonical isomorphisms between the spaces
of harmonic spinors with respect to a K\"ahler metric and certain
sheaf cohomology groups, as follows:
$$
ker(D) = H^{even} (X,\OO (\frac{1}{2}K)),
$$
$$
coker(D) = H^{odd} (X,\OO (\frac{1}{2}K)).
$$
This shows that the dimension of the space of harmonic spinors
is independent of the choice of K\"ahler metric compatible
with a fixed complex structure. The complex analytic version
of Problem 1 is:

\medskip\noindent
{\bf Problem 2.} {\sl If $X$ is a K\"ahlerian compact complex
manifold, does there exist a K\"ahler complex structure for
which there are no non--trivial harmonic spinors, i.e. such
that}
$$
\sum_i dim( H^i (\OO (\frac{1}{2}K)))=
\vert\chi(\OO (\frac{1}{2}K))\vert\ \ \ ?
$$
\medskip
\noindent
It is important to keep in mind that, unlike the cohomology
of $\OO (nK)$ for $n\in\Z$,
the cohomology of $\OO (\frac{1}{2}K)$ is not
invariant under deformations of the complex structure.
There are
(non--simply connected) counterexamples already in complex
dimension $1$.
However, the cohomology is semicontinuous, which implies that
if a complex structure exists with no non--trivial harmonic
spinors, then a generic complex structure also has this
property.

In the direction of Problem 2, Hitchin~\cite{hi} proved that a simply
connected spin K\"ahler surface which is not of general type
or is a complete intersection or a cyclic branched cover
of $\cp$ with smooth branch locus does not have non--trivial
harmonic spinors. Hitchin went on to conjecture that for
a generic complex structure on a simply connected
spin algebraic surface $H^1 (\OO (\frac{1}{2}K)) =0$. That this
cannot be true as stated follows from the existence of
simply connected spin algebraic surfaces for which
$\chi (\OO (\frac{1}{2}K)) <0$ due to Moishezon--Teicher~\cite{MT},
compare \cite{ko}. The correct
formulation of Hitchin's conjecture is that there should
be no non--trivial harmonic spinors for a generic complex
structure on a simply connected surface,
meaning $H^1 (\OO (\frac{1}{2}K)) =0$ if
$\chi (\OO (\frac{1}{2}K)) \geq 0$ and
$H^0 (\OO (\frac{1}{2}K)) = H^2 (\OO (\frac{1}{2}K)) =0$ if
$\chi (\OO (\frac{1}{2}K)) \leq 0$.

In \cite{ko}, we showed that the Moishezon--Teicher surfaces
have spaces of non--trivial harmonic spinors which can be
arbitrarily large. However, we do not know whether their
complex structures are generic.

We now show that Hitchin's conjecture is false, and that the
answer to Problem 2 is negative.

\begin{thm}
For each $k > 0$ there exist simply connected
 spin algebraic surfaces $X$ with
the following properties:

(i) the $dim H^i (\OO (\frac{1}{2}K))$ are deformation invariants
of $X$, and

(ii) $\sum_i dim(H^i (\OO (\frac{1}{2}K)))\geq
\vert\chi(\OO (\frac{1}{2}K))\vert +k$.
\end{thm}
\begin{pf}
Examples of such surfaces $X$ include the so--called
simple bihyperelliptic surfaces or bidouble covers
of $\CP^1\times\CP^1$ introduced by Catanese~\cite{ca}
and studied further by Catanese~\cite{ca2} and by
Manetti~\cite{ma}. A minimal
surface $X$ of general type is a bidouble cover
of $\CP^1\times\CP^1$ of type $(a,b)(m,n)$ if its
canonical model is defined in the total space of
$\OO (a,b)\oplus \OO (m,n)$ over $\CP^1\times\CP^1$
by a pair of equations
$$
z^2 = f(x,y)\ ,\ \ \ w^2 = g(x,y)\ ,
$$
where $f$ and $g$ are bihomogeneous polynomials
of bidegrees $(2a,2b)$ and $(2m,2n)$ respectively.

We consider the case where each of $a,b,m,n$ is at
least $4$. In this case Catanese~\cite{ca} showed that
the surfaces are simply connected and satisfy
\begin{equation}\label{K}
K_X =\pi^*\OO_{\CP^1\times\CP^1}(n+a-2,m+b-2)
\end{equation}
\begin{equation}\label{chi}
\chi (\OO_X)=2(n+a-1)(m+b-1)+2-\frac{1}{2}(m+b)(n+a)
+\frac{1}{2}(a-n)(b-m) .
\end{equation}
Moreover,
if $a\geq max\{ 2n+1,\ b+2\}$ and
$m\geq max\{ 2b+1,\ n+2\}$, then these surfaces form
an irreducible component \cite{ca2}, and even a
connected component \cite{ma} of the moduli space of surfaces
of general type. Thus, all the deformations of a bidouble
cover are bidouble covers of the same type.

If $a+n$ and $b+m$ are even, then $K_X$
is $2$--divisible and so $X$ is spin. As it is simply
connected, it has a unique spin structure, or half--canonical
divisor $\frac{1}{2}K_X$. From equations~\eqref{K} and
\eqref{chi} we obtain
\begin{alignat*}{1}
\chi (\OO (\frac{1}{2}K)) &= \chi (\OO_X)-\frac{1}{8}K^2\\
&= \frac{1}{2}(n+a)(m+b)+\frac{1}{2}(a-n)(b-m)\\
&=ab+mn> 0\ .
\end{alignat*}
On the other hand, we have
\begin{alignat*}{1}
dim H^0(\OO (\frac{1}{2}K)) &= dim H^0 (\OO_{\CP^1\times\CP^1}
(\frac{n+a-2}{2},\frac{m+b-2}{2}))\\
&=\frac{1}{4}(n+a)(m+b)\ .
\end{alignat*}
By Serre duality, we have $dim H^2(\OO (\frac{1}{2}K))=
dim H^0(\OO (\frac{1}{2}K))$. We conclude that
\begin{alignat*}{1}
\sum_i dim(H^i(\OO (\frac{1}{2}K)))-\vert\chi(\OO (\frac{1}{2}K))\vert
&= 2dim H^1(\OO (\frac{1}{2}K))\\ &=
2(-\chi (\OO (\frac{1}{2}K))+2dim(H^0 (\OO (\frac{1}{2}
K))))\\ &=
(a-n)(m-b)>0\ .
\end{alignat*}
By choosing $a,b,m,n$ appropriately, we can make this
arbitrarily large.
\end{pf}

\noindent
{\sl Acknowledgement:} The author wishes to thank F. Catanese
for some helpful comments.
\newpage
\bibliographystyle{amsplain}

\begin{thebibliography}{10}

\bibitem{ca} F.~Catanese, {\em On the moduli spaces of surfaces
of general type}, J. Diff. Geo. {\bf 19} (1984), 483--515.

\bibitem{ca2} F.~Catanese, {\em Automorphisms of rational
double points and moduli spaces of surfaces of general type},
Comp. Math. {\bf 61} (1987), 81--102.

\bibitem{hi} N.~Hitchin, {\em Harmonic Spinors}, Adv. in Math.
{\bf 14} (1974), 1--55.

\bibitem{ko} D.~Kotschick, {\em Non--trivial harmonic spinors on
certain algebraic surfaces}, in ``Einstein metrics and Yang--Mills
connections'', ed. T.~Mabuchi and S.~Mukai, Marcel Dekker, New York,
Basel, Hong Kong 1993.

\bibitem{ma} M.~Manetti, {\em On some components of moduli space
of surfaces of general type}, Comp. Math. {\bf 92} (1994), 285--297.

\bibitem{MT} B.~Moishezon and M.~Teicher, {\em Simply--connected
algebraic surfaces of positive index}, Invent. Math. {\bf 89}
(1987), 601--643.

\bibitem{wi} E.~Witten, {\em Monopoles and $4$--manifolds}, Math.
Research Letters {\bf 1} (1994), 769--796.

\end{thebibliography}

\end{document}